\documentclass[12pt]{article}
\usepackage{amssymb,amsmath,latexsym,textcomp,amsthm}
\usepackage{multirow}
\usepackage{rotating}
\usepackage{float}
\usepackage{pdflscape}
\usepackage{fancyhdr}
\usepackage{times}
\usepackage{graphicx} 
\graphicspath{ {images/} }
\begin{document}

\title{Finite Difference Weerakoon-Fernando Method to solve nonlinear equations without using derivatives}
\author{\small S.L. Heenatigala$^{1}$, Sunethra Weerakoon$^{2}$, T. G. I. Fernando$^{3}$ \\ \small Department of Mathematics$^{1,2}$ \& Department of Computer Science$^{3}$ \\ \small University of Sri Jayewardenepura, Gangodawila, Nugegoda, Sri Lanka}
\date{}
\maketitle

\begin{center}
   \textbf{Abstract }
\end{center}

\begin{sloppypar}This research was mainly conducted to explore the possibility of formulating an efficient algorithm to find roots of nonlinear equations without using the derivative of the function. The Weerakoon-Fernando method had been taken as the base in this project to find a new method without the derivative since Weerakoon-Fernando method gives 3rd order convergence.  After several unsuccessful attempts we were able to formulate the Finite Difference Weerakoon-Fernando Method (FDWFM) presented here. We noticed that the FDWFM approaches the root faster than any other existing method in the absence of the derivatives as an example, the popular nonlinear equation solver such as secant method (order of convergence is 1.618) in the absence of the derivative.  And the FDWFM had three function evaluations and  secant  method  had  two  function  evaluations. By implementing FDWFM on nonlinear equations with complex roots and also on systems of nonlinear equations, we received very encouraging results. When applying the FDWFM to systems of nonlinear equations, we resolved the involvement of the Jacobian problem by following the procedure in the Broyden\textquotesingle s method. The computational order of convergence of the FDWFM was close to 2.5 for all these cases. This will undoubtedly provide scientists the efficient numerical algorithm, that doesn't need the derivative of the function to solve nonlinear equations, that they were searching for over centuries. \\ 
\hspace{1.5cm}



\end{sloppypar}

\section{Introduction}
Even though the Newton\textquotesingle s method is popular in finding a single root of a nonlinear equation, the Weerakoon-Fernando method Fernando, (1998) and Weerakoon and Fernando, (2000) is better because its order of convergence is higher. Both the Newton\textquotesingle s method and the Weerakoon-Fernando method require the presence of the derivative of the function in the iterative equation. However, it so happens that some modelling situations do not provide value of the derivative baring the use of most computer algorithms using such data. The secant method is obtained by replacing the derivative of Newton\textquotesingle s method by a difference approximation. But the order of convergence of the secant method is 1.618 and the number of function evaluations required to carry out one iteration is two. However, this speed seems to be not adequate when considering the requirements of the current world with the high demand for more efficient algorithms due to the rapid development of technology. Thus we were compelled to search for a method to fulfill the current needs by replacing the derivatives in the fast converging Weerakoon-Fernando method by appropriate difference approximations due to the non-availability of powerful numerical methods that don't require derivatives to approximate roots of nonlinear equation. \\ This paper describes how to construct a best suitable method leading to an efficient algorithm without the derivatives to numerically solve nonlinear equations while optimizing the high order of convergence and the number of function evaluations. Then it describes how we checked that the method is in fact acceptable for nonlinear equations with complex roots. Finally it explains how that method was adapted even for systems of nonlinear equations, producing excellent results.\\ \\
\maketitle \textbf{\large 1.1 Preliminary Results}\\ \\
\maketitle \textbf{1.1.1 Order of convergence} \\ 
\maketitle \textbf{Definition 1.1.1 Order of Convergence of an iterative scheme
} \\ 
Let the iterative sequence: $ \left\{x_{n}: n = 1, 2,\ldots \right\} $ 
 that converges to $\underline{x}^{*}$ be generated by a numerical scheme. If there exists a constant $c\geq 0$, 
an integer $n_{0} \geq 0 $ and $\rho \geq 0$ such that for all $n > n_{0} $, 
the inequality below holds for any vector norm $\left\|.\right\| $.
\begin{align}
\left\| \underline{x} _{n+1} - \underline{x} ^{*} \right\| \leq c \left\| \underline{x}_{n}-\underline{x}^{*}\right\|^{\rho} {}\tag{2.1}
\end{align}Then the iterative scheme is said to be converge to $\underline{x}^{*}$ with $\rho^{th}$ order convergence. \\ \\
\maketitle \textbf{1.1.2 Computational Order of Convergence} \\ 
\textbf{Definition 1.1.2} \\  Let $\underline{x}^{*}$ be a root of the equation $\underline{F}(\underline{x}) = \underline{0} $ and suppose that $ \underline{x}_{n-1} $, $\underline{x}_{n} $ and $\underline{x}_{n+1}$ be consecutive iterates closer to the root $\underline{x}^{*}$, generated by an iterative scheme. Then the Computational Order of Convergence (COC) $\rho$ of the iterative scheme or the numerical algorithm can be approximated by

\begin{equation}
\rho \approx \frac{ln  \frac{\left\| \underline{x} _{n+1} - \underline{x} ^{*} \right\|}{\left\| \underline{x} _{n} - \underline{x} ^{*} \right\|} }{ln  \frac{\left\| \underline{x} _{n} - \underline{x} ^{*} \right\|}{\left\| \underline{x} _{n-1} - \underline{x} ^{*} \right\|} } \tag{2.2}
\end{equation}\\

\maketitle \textbf{\large 1.2 Numerical Schemes}\\

\maketitle \textbf{1.2.1 Weerakoon-Fernando Method to solve nonlinear equations in one variable}\\ 
Fernando, (1998) and Weerakoon and Fernando, (2000) introduced a third order convergent Weerakoon-Fernando Method (WFM) to solve nonlinear equations. The local model $ M_{n} \left(x\right) $ of WFM is given by the equation.

\begin{equation}
 M_{n}\left(x\right)=f\left(x_{n}\right)+\frac{1}{2}\left(x-x_{n}\right)\left[f^{'}\left(x_{n}\right)+f^{'}\left(x\right)\right] \tag{3.1}
\end{equation} When $x$ is taken as the next iterate $x_{n+1}$ and the root, we get the following formula

\begin{equation}
x_{n+1}=x_{n}-\frac{2f\left(x_{n}\right)}{f^{'}\left(x_{n}\right)+f^{'}\left(x_{n+1}\right)}\tag{3.2}   
\end{equation}here $n = 0, 1, 2,\ldots$ \\ \\
\maketitle\textbf{1.2.2 Weerakoon-Fernando formula for real variables}

\begin{equation}
x_{n+1}=x_{n}-\frac{2f\left(x_{n}\right)}{f^{'}\left(x_{n}\right)+f^{'}\left(x_{n+1}^{*}\right)}\tag{3.3}   
\end{equation} \\
\newpage
where   
\begin{equation}
 x_{n+1}^{*}=x_{n}-\frac{f\left(x_{n}\right)}{f^{'}\left(x_{n}\right)}\tag{3.4}  
\end{equation}here $n = 0, 1, 2,\ldots$ \\ \\

\maketitle\textbf{1.2.3 Weerakoon-Fernando formula for complex variables}\\ 
In Weerakoon-Fernando formula replace the variable $x$ with the complex variable $z$ in both sides then the following is the complex form of the WFM.
\begin{equation}
z_{n+1}=z_{n}-\frac{2f\left(z_{n}\right)}{f^{'}\left(z_{n}\right)+f^{'}\left(z_{n+1}^{*}\right)}\tag{3.5}   
\end{equation} where   
\begin{equation}
 z_{n+1}^{*}=z_{n}-\frac{f\left(z_{n}\right)}{f^{'}\left(z_{n}\right)}\tag{3.6}  
\end{equation}here $n = 0, 1, 2,\ldots$ \\ \\
\maketitle\textbf{1.2.4 The Weerakoon-Fernando Method to solve systems of nonlinear equations}\\ 
Nishani, (2015) and Nishani, Weerakoon, Fernando and Liyanage, (2018) introduced a third order convergent Weerakoon-Fernando Method (WFM) to solve systems of nonlinear equations. When $\underline{F}$ is a vector valued function with non-zero derivatives defined on the set $D\subset \Re ^{n} $ and $x$, $x_{0} \in D$, the extension of the WFM to solve systems of nonlinear equations can be given as follows. 

\begin{equation}
\underline{x}_{n+1}=\underline{x}_{n}-2\left[J\left(\underline{F}\left(\underline{x}_{n}\right)\right)+J\left(\underline{F}\left(\underline{x}_{n+1}^{\lambda}\right)\right)\right]^{-1}\left(\underline{F}\left(\underline{x}_{n}\right)\right)\tag{3.7}
\end{equation}where\;

\begin{equation}
\underline{x}_{n+1}^{\lambda}=\underline{x}_{n}-\left[J\left(\underline{F}\left(\underline{x}_{n}\right)\right)\right]^{-1}\left(\underline{F}\left(\underline{x}_{n}\right)\right)\tag{3.8}
\end{equation}
here $n = 0, 1, 2,\ldots$ and $\underline{x}^{n}$ is the $n^{th}$ iterate and $J$ is the Jacobian matrix of $\underline{F}$. \\ \\
\newpage
 \maketitle\textbf{1.2.5 Secant Method for Systems of Nonlinear Equations} \\ 
When we use the secant method for systems of nonlinear equations we face a special problem with the Jacobian. In secant method for systems, a vector is present in the denominator. To overcome the problem of having to take the inverse of a vector we follow the most popular secant approximation proposed by C. Broyden The algorithm is analogous to Newton\textquotesingle s method, but it replaces the analytic Jacobian by the following approximation. This method is called the Broyden\textquotesingle s method Atkinson, (1988). \\ \\ 

\maketitle ALGORITHEM  \textbf{Broyden\textquotesingle s method}. \\ \\
Given $\underline{F}:\Re^{n}\rightarrow \Re ^{n}, \underline{x}_{0}\in \Re ^{n}, A^{0}\in \Re ^{n\times n}$ \\ $s_{k}$- Initial Step, $y_{k}$- Yield of current Step, $\underline{x}_{k+1}$- Next iteration  \\ FOR $k=0$ to $a$ \\ $s_{k}=x_{k+1}-x_{k}$ \\ ${y}_{k}=\underline{F}\left(\underline{x}_{k+1}\right)-\underline{F}\left(\underline{x}_{k}\right) $ for $A_{k}s_{k}={y}_{k} $ 
\\$A_{k}\underline{s}_{k}=-\underline{F}\left(\underline{x}_{k}\right)$ \\ $\underline{x}_{k+1}=\underline{x}_{k}+\underline{s}_{k}$  \\ $\underline{y}_{k}=\underline{F}\left(\underline{x}_{k+1}\right)-\underline{F}\left(\underline{x}_{k}\right)$ \\ $A_{k+1}=A_{k}+\frac{\underline{y}_{k}-A_{k}\underline{s}_{k}}{\underline{s}_{k}^{t}\underline{s}_{k}}\underline{s}_{k}^{t}$ \\ END FOR \\ Here the final step is used to replace the analytic Jacobian by a matrix. \\

\section{Methodology}
\maketitle \textbf{ 2.1 Derivation of FDWFM from WFM}\\ \\
Here we use both backward and forward difference approximations for the derivatives in the same formula appropriately to get an acceptable result. \\
Consider the Weerakoon-Fernando Method 
\begin{equation}
x_{n+1}=x_{n}-\frac{2f\left(x_{n}\right)}{f^{'}\left(x_{n}\right)+f^{'}\left(x_{n+1}\right)} \tag{4.1}
\end{equation}Substituting the forward difference approximation for

\begin{equation}
f^{'}\left(x_{n}\right)\approx\frac{f\left(x_{n+1}\right)-f\left(x_{n}\right)}{x_{n+1}-x_{n}} \tag{4.2}
\end{equation}and the backward difference approximation for 

\begin{equation}
f^{'}\left(x_{n+1}\right)\approx \frac{f\left(x_{n+1}\right)-f\left(x_{n}\right)}{x_{n+1}-x_{n}} \tag{4.3}
\end{equation}in equation $\left(4.1\right)$, we get
\begin{equation}
x_{n+1}=x_{n}-\frac{2f\left(x_{n}\right)}{\frac{f\left(x_{n+1}\right)-f\left(x_{n}\right)}{x_{n+1}-x_{n}} +\frac{f\left(x_{n+1}\right)-f\left(x_{n}\right)}{x_{n+1}-x_{n}} }\tag{4.4}
\end{equation}

\begin{equation}
x_{n+1}=x_{n}-\frac{2f\left(x_{n}\right)\left(x_{n+1}-x_{n}\right)}{2\left(f\left(x_{n+1}\right)-f\left(x_{n}\right)\right)}\tag{4.5}
\end{equation}Thus the new iterative formula without derivative terms is
\begin{equation}
x_{n+1}=x_{n}-\frac{f\left(x_{n}\right)\left(x_{n+1}-x_{n}\right)}{f\left(x_{n+1}\right)-f\left(x_{n}\right)}\tag{4.6}
\end{equation}But it is an implicit method requiring $x_{n+1}$ term at the $\left(n+1\right)^{th}$  iterative step to calculate the $\left(n+1\right)^{th}$ iterate itself. The secant method can be used to replace the term $x_{n+1}$ on the RHS of the above equation to overcome this difficulty. Thus the Finite Difference Weerakoon-Fernando Method (FDWFM) is obtained as
\begin{equation}
x_{n+1}=x_{n}-\frac{f\left(x_{n}\right)\left(x_{n+1}^{*}-x_{n}\right)}{f\left(x_{n+1}^{*}\right)-f\left(x_{n}\right)}\tag{4.7}
\end{equation}where
\begin{equation}
x_{n+1}^{*}=x_{n}-\frac{f\left(x_{n}\right)\left(x_{n}-x_{n-1}\right)}{f\left(x_{n}\right)-f\left(x_{n-1}\right)}\tag{4.8}
\end{equation}$n=1,2,\ldots$ \\ \\  \\ \\ \\ \\
\maketitle \textbf{2.2 Application of FDWFM for nonlinear equations with complex roots} \\ 
In the Finite Difference Weerakoon-Fernando formula replace the variable x with the complex variable z in both sides to get the following form: \\
\begin{equation}
z_{n+1}=z_{n}-\frac{f\left(z_{n}\right)\left(z_{n+1}^{*}-z_{n}\right)}{f\left(z_{n+1}^{*}\right)-f\left(z_{n}\right)}\tag{4.9}
\end{equation}where
\begin{equation}
z_{n+1}^{*}=z_{n}-\frac{f\left(z_{n}\right)\left(z_{n}-z_{n-1}\right)}{f\left(z_{n}\right)-f\left(z_{n-1}\right)}\tag{4.10}
\end{equation}
Here $n =0,1,2,\ldots$ \\ \\ 

\maketitle\textbf{2.3 Application of FDWFM for systems of nonlinear equations} \\   
Choose an initial estimate $ \underline  {x}_{0} \in \Re ^{n}$ and a non-singular initial Matrix $A^{0} \in \Re ^{n\times n}$. Set $k:=0$ and repeat the following sequence of steps until $\left\|\underline{F}(\underline{x}^{k})\right\| < tolerance $ 
\begin{enumerate}
	\item Solve $A_{k}\underline{I}_{k}=\underline{F}\left(\underline{X}^{*}_{k+1}\right)-\underline{F}\left(\underline{X}_{k}\right) \ \ \ for \ \underline{I}_{k} \ \ \& \ \ \underline{X}^{*}_{k+1}$ from Broyden\textquotesingle s method 
	\item $\underline{X}_{k+1}=\underline{X}_{k}+\underline{I}_{k}$
	\item $\overline{\underline{Y}}_{k}=\underline{F}\left(\underline{X}_{k+1}\right)-\underline{F}\left(\underline{X}_{k}\right)$ 
	\item $A_{k+1}=A_{k}+\frac{\left(\overline{\underline{Y}}_{k}-A_{k}\underline{I}_{k}\right)\underline{I}_{k}^{t}}{\underline{I}_{k}^{t}\underline{I}_{k}}$
\end{enumerate}

\section{Results \& Discussion}
Results are given in Table 1, Table 2, Table 3, Table 4, Table 5, Table 6 and Table 7.

\begin{sidewaystable}[]
\centering
\caption{Comparison of FDWFM with secant method \& Newton\textquotesingle s for uni-variate nonlinear equations with real roots}
\label{my-label}
\begin{tabular}{ccccccccccccc}
\hline
\multicolumn{1}{|c|}{\multirow{2}{*}{Function}} & \multicolumn{1}{l|}{\multirow{2}{*}{X0}} & \multicolumn{1}{l|}{\multirow{2}{*}{X1}} & \multicolumn{3}{c|}{i} & \multicolumn{3}{c|}{COC} & \multicolumn{3}{c|}{NFC} & \multicolumn{1}{c|}{\multirow{2}{*}{Root}} \\ \cline{4-12}
\multicolumn{1}{|c|}{} & \multicolumn{1}{l|}{} & \multicolumn{1}{l|}{} & \multicolumn{1}{c|}{SCT} & \multicolumn{1}{c|}{NM} & \multicolumn{1}{c|}{FDWFM} & \multicolumn{1}{c|}{SCT} & \multicolumn{1}{c|}{NM} & \multicolumn{1}{c|}{FDWFM} & \multicolumn{1}{c|}{SCT} & \multicolumn{1}{c|}{NM} & \multicolumn{1}{c|}{FDWFM} & \multicolumn{1}{c|}{} \\ \hline
\multicolumn{1}{|c|}{$x^{3}+5x+4$} & \multicolumn{1}{c|}{0} & \multicolumn{1}{c|}{1} & \multicolumn{1}{c|}{8} & \multicolumn{1}{c|}{6} & \multicolumn{1}{c|}{4} & \multicolumn{1}{c|}{1.365} & \multicolumn{1}{c|}{1.97} & \multicolumn{1}{c|}{2.48791} & \multicolumn{1}{c|}{9} & \multicolumn{1}{c|}{8} & \multicolumn{1}{c|}{9} & \multicolumn{1}{c|}{0.7240755} \\ \hline
\multicolumn{1}{|c|}{$4cosx+e^{x}$} & \multicolumn{1}{c|}{0} & \multicolumn{1}{c|}{1.5} & \multicolumn{1}{c|}{10} & \multicolumn{1}{c|}{7} & \multicolumn{1}{c|}{4} & \multicolumn{1}{c|}{1.904} & \multicolumn{1}{c|}{1.99} & \multicolumn{1}{c|}{2.39479} & \multicolumn{1}{c|}{11} & \multicolumn{1}{c|}{9} & \multicolumn{1}{c|}{9} & \multicolumn{1}{c|}{0.904788} \\ \hline
\multicolumn{1}{|c|}{$sin^{2}x-x^{2}+1$} & \multicolumn{1}{c|}{1} & \multicolumn{1}{c|}{3} & \multicolumn{1}{c|}{9} & \multicolumn{1}{c|}{6} & \multicolumn{1}{c|}{5} & \multicolumn{1}{c|}{1.404} & \multicolumn{1}{c|}{1.98} & \multicolumn{1}{c|}{2.39405} & \multicolumn{1}{c|}{10} & \multicolumn{1}{c|}{8} & \multicolumn{1}{c|}{11} & \multicolumn{1}{c|}{1.404491} \\ \hline
\multicolumn{1}{|c|}{$x^{2}-e^{x}-3x+2$} & \multicolumn{1}{c|}{1} & \multicolumn{1}{c|}{2} & \multicolumn{1}{c|}{8} & \multicolumn{1}{c|}{6} & \multicolumn{1}{c|}{4} & \multicolumn{1}{c|}{1.257} & \multicolumn{1}{c|}{1.56} & \multicolumn{1}{c|}{2.48162} & \multicolumn{1}{c|}{9} & \multicolumn{1}{c|}{9} & \multicolumn{1}{c|}{9} & \multicolumn{1}{c|}{0.2575302} \\ \hline
\multicolumn{1}{|c|}{$cosx-x$} & \multicolumn{1}{c|}{0} & \multicolumn{1}{c|}{1} & \multicolumn{1}{c|}{6} & \multicolumn{1}{c|}{5} & \multicolumn{1}{c|}{3} & \multicolumn{1}{c|}{1.739} & \multicolumn{1}{c|}{1.99} & \multicolumn{1}{c|}{Not Defined} & \multicolumn{1}{c|}{7} & \multicolumn{1}{c|}{7} & \multicolumn{1}{c|}{7} & \multicolumn{1}{c|}{0.739085} \\ \hline
\multicolumn{1}{|c|}{$\left(1-x\right)^{3}-1$} & \multicolumn{1}{c|}{2.5} & \multicolumn{1}{c|}{3.5} & \multicolumn{1}{c|}{9} & \multicolumn{1}{c|}{8} & \multicolumn{1}{c|}{5} & \multicolumn{1}{c|}{1.864} & \multicolumn{1}{c|}{1.98} & \multicolumn{1}{c|}{2.37182} & \multicolumn{1}{c|}{10} & \multicolumn{1}{c|}{10} & \multicolumn{1}{c|}{10} & \multicolumn{1}{c|}{2.00000} \\ \hline
\multicolumn{1}{|c|}{$xe^{x^{2}}-sin^{2}x+3cosx+5$} & \multicolumn{1}{c|}{2} & \multicolumn{1}{c|}{2.5} & \multicolumn{1}{c|}{9} & \multicolumn{1}{c|}{7} & \multicolumn{1}{c|}{4} & \multicolumn{1}{c|}{1.562} & \multicolumn{1}{c|}{1.99} & \multicolumn{1}{c|}{2.53667} & \multicolumn{1}{c|}{10} & \multicolumn{1}{c|}{9} & \multicolumn{1}{c|}{10} & \multicolumn{1}{c|}{2.154434} \\ \hline
\multicolumn{1}{|c|}{$x^{2}sin^{2}x+e^{x^{2}cosxsinx}-28$} & \multicolumn{1}{c|}{-2} & \multicolumn{1}{c|}{-1} & \multicolumn{1}{c|}{11} & \multicolumn{1}{c|}{9} & \multicolumn{1}{c|}{5} & \multicolumn{1}{c|}{1.432} & \multicolumn{1}{c|}{1.99} & \multicolumn{1}{c|}{2.49479} & \multicolumn{1}{c|}{12} & \multicolumn{1}{c|}{11} & \multicolumn{1}{c|}{12} & \multicolumn{1}{c|}{-1.207647} \\ \hline
\multicolumn{1}{|c|}{$e^{x^{2}+7x-30}-1$} & \multicolumn{1}{c|}{4} & \multicolumn{1}{c|}{5} & \multicolumn{1}{c|}{29} & \multicolumn{1}{c|}{25} & \multicolumn{1}{c|}{14} & \multicolumn{1}{c|}{1.781} & \multicolumn{1}{c|}{1.99} & \multicolumn{1}{c|}{2.41052} & \multicolumn{1}{c|}{30} & \multicolumn{1}{c|}{30} & \multicolumn{1}{c|}{30} & \multicolumn{1}{c|}{3.4374717} \\ \hline
\multicolumn{1}{|c|}{$x^{3}-10$} & \multicolumn{1}{c|}{2.9} & \multicolumn{1}{c|}{3.5} & \multicolumn{1}{c|}{12} & \multicolumn{1}{c|}{11} & \multicolumn{1}{c|}{6} & \multicolumn{1}{c|}{1.653} & \multicolumn{1}{c|}{1.99} & \multicolumn{1}{c|}{2.41272} & \multicolumn{1}{c|}{13} & \multicolumn{1}{c|}{13} & \multicolumn{1}{c|}{13} & \multicolumn{1}{c|}{3.000000} \\ \hline
\multicolumn{6}{l}{FDWFM-Finite Difference Weerakoon-Fernando Method} & \multicolumn{7}{l}{SCT - Secant method} \\
\multicolumn{6}{l}{NM-Newton\textquotesingle s method} & \multicolumn{7}{l}{i-Number of iterations to approximate the root} \\
\multicolumn{6}{l}{COC-Computational Order of Convergence} & \multicolumn{7}{l}{NFE-Number of Function Evaluations}
\end{tabular}
\end{sidewaystable}

\begin{sidewaystable}[]
\centering
\caption{Comparison of FDWFM with secant method \& Newton\textquotesingle s method for nonlinear equations with complex roots}
\label{my-label}
\begin{tabular}{ccccccccccccc}
\hline
\multicolumn{1}{|c|}{\multirow{2}{*}{Function}} & \multicolumn{1}{c|}{\multirow{2}{*}{X0}} & \multicolumn{1}{c|}{\multirow{2}{*}{X1}} & \multicolumn{3}{c|}{i} & \multicolumn{3}{c|}{COC} & \multicolumn{3}{c|}{NFC} & \multicolumn{1}{c|}{\multirow{2}{*}{Root}} \\ \cline{4-12}
\multicolumn{1}{|c|}{} & \multicolumn{1}{c|}{} & \multicolumn{1}{c|}{} & \multicolumn{1}{c|}{SCT} & \multicolumn{1}{c|}{NM} & \multicolumn{1}{c|}{\begin{tabular}[c]{@{}c@{}}FDW\\ FM\end{tabular}} & \multicolumn{1}{c|}{SCT} & \multicolumn{1}{c|}{NM} & \multicolumn{1}{c|}{\begin{tabular}[c]{@{}c@{}}FDW\\ FM\end{tabular}} & \multicolumn{1}{c|}{SCT} & \multicolumn{1}{c|}{NM} & \multicolumn{1}{c|}{\begin{tabular}[c]{@{}c@{}}FDW\\ FM\end{tabular}} & \multicolumn{1}{c|}{} \\ \hline
\multicolumn{1}{|c|}{$z^{2}+1$} & \multicolumn{1}{c|}{0,0.5} & \multicolumn{1}{c|}{0.1,0.8} & \multicolumn{1}{c|}{8} & \multicolumn{1}{c|}{7} & \multicolumn{1}{c|}{4} & \multicolumn{1}{c|}{1.565} & \multicolumn{1}{c|}{2.008} & \multicolumn{1}{c|}{2.43010} & \multicolumn{1}{c|}{9} & \multicolumn{1}{c|}{9} & \multicolumn{1}{c|}{9} & \multicolumn{1}{c|}{0 + 1.000i} \\ \hline
\multicolumn{1}{|c|}{$z^{2}+e^{z}-3z-3$} & \multicolumn{1}{c|}{1,0.5} & \multicolumn{1}{c|}{1.2,0.7} & \multicolumn{1}{c|}{7} & \multicolumn{1}{c|}{6} & \multicolumn{1}{c|}{4} & \multicolumn{1}{c|}{1.704} & \multicolumn{1}{c|}{1.996} & \multicolumn{1}{c|}{2.40123} & \multicolumn{1}{c|}{8} & \multicolumn{1}{c|}{9} & \multicolumn{1}{c|}{9} & \multicolumn{1}{c|}{-0.8800+0.0000i} \\ \hline
\multicolumn{1}{|c|}{$sin^{2}z-z-2$} & \multicolumn{1}{c|}{1,0.5} & \multicolumn{1}{c|}{1.2,0.7} & \multicolumn{1}{c|}{12} & \multicolumn{1}{c|}{10} & \multicolumn{1}{c|}{6} & \multicolumn{1}{c|}{1.504} & \multicolumn{1}{c|}{2.001} & \multicolumn{1}{c|}{2.36646} & \multicolumn{1}{c|}{13} & \multicolumn{1}{c|}{12} & \multicolumn{1}{c|}{13} & \multicolumn{1}{c|}{1.0037+ 1.1767i} \\ \hline
\multicolumn{1}{|c|}{$ze^{z^{2}}-sinz+3cosz+1$} & \multicolumn{1}{c|}{-2,1} & \multicolumn{1}{c|}{-1,1.5} & \multicolumn{1}{c|}{16} & \multicolumn{1}{c|}{13} & \multicolumn{1}{c|}{7} & \multicolumn{1}{c|}{1.657} & \multicolumn{1}{c|}{1.991} & \multicolumn{1}{c|}{2.42055} & \multicolumn{1}{c|}{17} & \multicolumn{1}{c|}{15} & \multicolumn{1}{c|}{15} & \multicolumn{1}{c|}{1.1199- 0.7028i} \\ \hline
\multicolumn{1}{|c|}{$z^{3}-3z^{2}+3z$} & \multicolumn{1}{c|}{1.5,0.5} & \multicolumn{1}{c|}{1.7,0.19} & \multicolumn{1}{c|}{8} & \multicolumn{1}{c|}{7} & \multicolumn{1}{c|}{4} & \multicolumn{1}{c|}{1.739} & \multicolumn{1}{c|}{1.996} & \multicolumn{1}{c|}{2.37201} & \multicolumn{1}{c|}{9} & \multicolumn{1}{c|}{9} & \multicolumn{1}{c|}{9} & \multicolumn{1}{c|}{1.6299+ 1.091i} \\ \hline
\multicolumn{1}{|c|}{$\left(1-z\right)^{3}+1$} & \multicolumn{1}{c|}{1.5,0.5} & \multicolumn{1}{c|}{1.5,1} & \multicolumn{1}{c|}{10} & \multicolumn{1}{c|}{9} & \multicolumn{1}{c|}{5} & \multicolumn{1}{c|}{1.664} & \multicolumn{1}{c|}{2.002} & \multicolumn{1}{c|}{2.40747} & \multicolumn{1}{c|}{11} & \multicolumn{1}{c|}{11} & \multicolumn{1}{c|}{11} & \multicolumn{1}{c|}{1.500 + 0.8660i} \\ \hline
\multicolumn{1}{|c|}{$z^{5}-z^{4}+7z^{3}-5z^{2}+4z-4$} & \multicolumn{1}{c|}{0,0.4} & \multicolumn{1}{c|}{0.1,0.5} & \multicolumn{1}{c|}{8} & \multicolumn{1}{c|}{7} & \multicolumn{1}{c|}{5} & \multicolumn{1}{c|}{1.562} & \multicolumn{1}{c|}{2.001} & \multicolumn{1}{c|}{2.20234} & \multicolumn{1}{c|}{9} & \multicolumn{1}{c|}{10} & \multicolumn{1}{c|}{11} & \multicolumn{1}{c|}{-0.0880 + 0.869i} \\ \hline
\multicolumn{1}{|c|}{$z^{4}+1$} & \multicolumn{1}{c|}{0.01,0.5} & \multicolumn{1}{c|}{0.3,0.8} & \multicolumn{1}{c|}{16} & \multicolumn{1}{c|}{15} & \multicolumn{1}{c|}{8} & \multicolumn{1}{c|}{1.532} & \multicolumn{1}{c|}{2.001} & \multicolumn{1}{c|}{2.43348} & \multicolumn{1}{c|}{17} & \multicolumn{1}{c|}{17} & \multicolumn{1}{c|}{17} & \multicolumn{1}{c|}{0.7071 + 0.7071i} \\ \hline
\multicolumn{1}{|c|}{$z^{2}sin^{2}2z+e^{z^{2}coszsinz}+10$} & \multicolumn{1}{c|}{2.5,0.3} & \multicolumn{1}{c|}{3,0.4} & \multicolumn{1}{c|}{16} & \multicolumn{1}{c|}{14} & \multicolumn{1}{c|}{7} & \multicolumn{1}{c|}{1.681} & \multicolumn{1}{c|}{2.003} & \multicolumn{1}{c|}{2.42235} & \multicolumn{1}{c|}{17} & \multicolumn{1}{c|}{16} & \multicolumn{1}{c|}{15} & \multicolumn{1}{c|}{3.1981 + 0.4354i} \\ \hline
\multicolumn{6}{l}{FDWFM-Finite Difference Weerakoon-Fernando Method} & \multicolumn{7}{l}{SCT - Secant method} \\
\multicolumn{6}{l}{NM-Newton\textquotesingle s method} & \multicolumn{7}{l}{i-Number of iterations to approximate the root} \\
\multicolumn{6}{l}{COC-Computational Order of Convergence} & \multicolumn{7}{l}{NFE-Number of Function Evaluations}
\end{tabular}
\end{sidewaystable}

\begin{sidewaystable}[]
\centering
\caption{Comparison of FDWFM with Broyden\textquotesingle s method for systems of two nonlinear equations}
\label{my-label}
\begin{tabular}{lcccccll}
\hline
\multicolumn{1}{|l|}{\multirow{2}{*}{}} & \multicolumn{1}{c|}{\multirow{2}{*}{Function}} & \multicolumn{1}{c|}{\multirow{2}{*}{\underline{X}0}} & \multicolumn{2}{c|}{i} & \multicolumn{2}{c|}{COC} & \multicolumn{1}{c|}{\multirow{2}{*}{Root}} \\ \cline{4-7}
\multicolumn{1}{|l|}{} & \multicolumn{1}{c|}{} & \multicolumn{1}{c|}{} & \multicolumn{1}{c|}{BMFSM} & \multicolumn{1}{c|}{FDWFM} & \multicolumn{1}{c|}{BMFSM} & \multicolumn{1}{c|}{FDWFM} & \multicolumn{1}{c|}{} \\ \hline
\multicolumn{1}{|l|}{1} & \multicolumn{1}{c|}{\begin{tabular}[c]{@{}c@{}}$x+y-3$\\ $x^{2}+y^{2}-9$\end{tabular}} & \multicolumn{1}{c|}{(1,3)} & \multicolumn{1}{c|}{12} & \multicolumn{1}{c|}{10} & \multicolumn{1}{c|}{1.678007} & \multicolumn{1}{l|}{2.6416} & \multicolumn{1}{l|}{\begin{tabular}[c]{@{}l@{}}0.00000064\\   2.99998136\end{tabular}} \\ \hline
\multicolumn{1}{|l|}{2} & \multicolumn{1}{c|}{\begin{tabular}[c]{@{}c@{}}$x^{4}+y^{4}-67$\\ $x^{3}-3xy^{2}+35$\end{tabular}} & \multicolumn{1}{c|}{(2,3)} & \multicolumn{1}{c|}{13} & \multicolumn{1}{c|}{12} & \multicolumn{1}{c|}{1.712084} & \multicolumn{1}{l|}{2.5687} & \multicolumn{1}{l|}{\begin{tabular}[c]{@{}l@{}}1.88364056\\   2.71594745\end{tabular}} \\ \hline
\multicolumn{1}{|l|}{3} & \multicolumn{1}{c|}{\begin{tabular}[c]{@{}c@{}}$x^{2}-10x+y^{2}+8$\\ $xy^{2}+x-10y+8$\end{tabular}} & \multicolumn{1}{c|}{(0,0)} & \multicolumn{1}{c|}{11} & \multicolumn{1}{c|}{8} & \multicolumn{1}{c|}{1.687346} & \multicolumn{1}{l|}{2.8984} & \multicolumn{1}{l|}{\begin{tabular}[c]{@{}l@{}}0.99999927\\   1.00000048\end{tabular}} \\ \hline
\multicolumn{1}{|l|}{4} & \multicolumn{1}{c|}{\begin{tabular}[c]{@{}c@{}}$x-cosy$\\ $sinx+0.5y$\end{tabular}} & \multicolumn{1}{c|}{(0,-0.5)} & \multicolumn{1}{c|}{11} & \multicolumn{1}{c|}{7} & \multicolumn{1}{c|}{1.78364} & \multicolumn{1}{l|}{2.2331} & \multicolumn{1}{l|}{\begin{tabular}[c]{@{}l@{}}0.53038868\\   -1.01173734\end{tabular}} \\ \hline
\multicolumn{1}{|l|}{5} & \multicolumn{1}{c|}{\begin{tabular}[c]{@{}c@{}}$x^{2}+y^{2}-2$\\ $e^{x-1}+y^{3}-2$\end{tabular}} & \multicolumn{1}{c|}{(0.5,0.5)} & \multicolumn{1}{c|}{15} & \multicolumn{1}{c|}{10} & \multicolumn{1}{c|}{1.69356} & \multicolumn{1}{l|}{2.2420} & \multicolumn{1}{l|}{\begin{tabular}[c]{@{}l@{}}0.99998777\\   1.00001497\end{tabular}} \\ \hline
\multicolumn{1}{|l|}{6} & \multicolumn{1}{c|}{\begin{tabular}[c]{@{}c@{}}$-x^{2}-x+2y-18$\\ $\left(x-1\right)^{2}+\left(y-6\right)2-25$\end{tabular}} & \multicolumn{1}{c|}{(-5,-5)} & \multicolumn{1}{c|}{24} & \multicolumn{1}{c|}{15} & \multicolumn{1}{c|}{1.537854} & \multicolumn{1}{l|}{2.3251} & \multicolumn{1}{l|}{\begin{tabular}[c]{@{}l@{}}1.54694636\\   10.9699948\end{tabular}} \\ \hline
\multicolumn{1}{|l|}{7} & \multicolumn{1}{c|}{\begin{tabular}[c]{@{}c@{}}$2cosy+7sinx-10x$\\ $7cosx-2siny-10y$\end{tabular}} & \multicolumn{1}{c|}{(0,0)} & \multicolumn{1}{c|}{10} & \multicolumn{1}{c|}{9} & \multicolumn{1}{c|}{1.62944} & \multicolumn{1}{l|}{2.4499} & \multicolumn{1}{l|}{\begin{tabular}[c]{@{}l@{}}0.52651702\\   0.50792810\end{tabular}} \\ \hline
\multicolumn{3}{l}{FDWFM - Finite Difference Weerakoon-Fernando Method} & \multicolumn{5}{l}{BMFSM - Broyden\textquotesingle s method for secant method} \\
\multicolumn{3}{l}{i - Number of iterations to approximate the root} & \multicolumn{5}{l}{COC - Computational Order of Convergence} \\
\end{tabular}
\end{sidewaystable}

\begin{sidewaystable}[]
\centering
\caption{Comparison of FDWFM with Broyden\textquotesingle s method for systems of three nonlinear equations}
\label{my-label}
\begin{tabular}{lcccccll}
\hline
\multicolumn{1}{|l|}{\multirow{2}{*}{}} & \multicolumn{1}{c|}{\multirow{2}{*}{Function}} & \multicolumn{1}{c|}{\multirow{2}{*}{\underline{X}0}} & \multicolumn{2}{c|}{i} & \multicolumn{2}{c|}{COC} & \multicolumn{1}{c|}{\multirow{2}{*}{Root}} \\ \cline{4-7}
\multicolumn{1}{|l|}{} & \multicolumn{1}{c|}{} & \multicolumn{1}{c|}{} & \multicolumn{1}{c|}{BMFSM} & \multicolumn{1}{c|}{FDWFM} & \multicolumn{1}{c|}{BMFSM} & \multicolumn{1}{c|}{FDWFM} & \multicolumn{1}{c|}{} \\ \hline
\multicolumn{1}{|l|}{1} & \multicolumn{1}{c|}{\begin{tabular}[c]{@{}c@{}}$x^{2}+y^{2}+z^{2}$\\ $x^{2}-y^{2}+z^{2}$\\ $x^{2}+y^{2}-z^{2}$\end{tabular}} & \multicolumn{1}{c|}{(1,1,1)} & \multicolumn{1}{c|}{15} & \multicolumn{1}{c|}{14} & \multicolumn{1}{c|}{1.71123} & \multicolumn{1}{l|}{2.3883} & \multicolumn{1}{l|}{\begin{tabular}[c]{@{}l@{}}0.00270271\\   0.00270271\\   0.00270271\end{tabular}} \\ \hline
\multicolumn{1}{|l|}{2} & \multicolumn{1}{c|}{\begin{tabular}[c]{@{}c@{}}$3x^{2}cos\left(?yz\right)-1/2$\\ $x^{2}-81\left(y+0.1\right)^{2}+sin?\left(z\right)+1.06$\\ $e^{-\left(xy\right)}+20z+\left(10\pi-3\right)/3$\end{tabular}} & \multicolumn{1}{c|}{(1,1,-1)} & \multicolumn{1}{c|}{16} & \multicolumn{1}{c|}{13} & \multicolumn{1}{c|}{1.84657} & \multicolumn{1}{l|}{2.81152} & \multicolumn{1}{l|}{\begin{tabular}[c]{@{}l@{}}0.70712121\\   0.01416426\\   -0.52299771\end{tabular}} \\ \hline
\multicolumn{1}{|l|}{3} & \multicolumn{1}{c|}{\begin{tabular}[c]{@{}c@{}}$x+e^{x-1}+\left(y+z\right)^{2}-27$\\ $ e^{y-2}/x+z^{2}-10$\\ $y^{2}+sin\left(?y-2\right)+z-7$\end{tabular}} & \multicolumn{1}{c|}{(1.4,2.2,3.1)} & \multicolumn{1}{c|}{9} & \multicolumn{1}{c|}{8} & \multicolumn{1}{c|}{1.877212} & \multicolumn{1}{l|}{2.79400} & \multicolumn{1}{l|}{\begin{tabular}[c]{@{}l@{}}0.99996491\\   2.00009396\\   2.99980941\end{tabular}} \\ \hline
\multicolumn{1}{|l|}{4} & \multicolumn{1}{c|}{\begin{tabular}[c]{@{}c@{}}$15x+y^{2}-4z-13$\\ $x^{2}+10y-z-11$\\ $y^{3}-25z+22$\end{tabular}} & \multicolumn{1}{c|}{(3,3,2)} & \multicolumn{1}{c|}{14} & \multicolumn{1}{c|}{12} & \multicolumn{1}{c|}{1.722963} & \multicolumn{1}{l|}{2.94914} & \multicolumn{1}{l|}{\begin{tabular}[c]{@{}l@{}}1.03640452\\   1.08570343\\   0.93119446\end{tabular}} \\ \hline
\multicolumn{3}{l}{FDWFM - Finite Difference Weerakoon-Fernando Method} & \multicolumn{5}{l}{BMFSM - Broyden\textquotesingle s method for secant method} \\
\multicolumn{3}{l}{i - Number of iterations to approximate the root} & \multicolumn{5}{l}{COC - Computational Order of Convergence} \\
\end{tabular}
\end{sidewaystable}

\begin{sidewaystable}[]
\centering
\caption{Comparison of FDWFM with Newton\textquotesingle s method and Broyden\textquotesingle s method for systems of four nonlinear equations $f_{1}=x_{1}+x_{2}-2, f_{2}=x_{1}x_{3}+x_{2}x_{4}, f_{3}=x_{1}x_{3}^{2}+x_{2}x_{4}^{2}-2/3, f_{4}=x_{1}x_{3}^{3}+x_{2}x_{4}^{3}$}
\label{my-label}
\begin{tabular}{cclccccc}
\hline
\multicolumn{1}{|c|}{\multirow{2}{*}{\begin{tabular}[c]{@{}c@{}}Initial guess\\ X0\end{tabular}}} & \multicolumn{3}{c|}{i} & \multicolumn{3}{c|}{COC} & \multicolumn{1}{c|}{\multirow{2}{*}{Root}} \\ \cline{2-7}
\multicolumn{1}{|c|}{} & \multicolumn{1}{c|}{BMFSM} & \multicolumn{1}{l|}{NM} & \multicolumn{1}{c|}{FDWFM} & \multicolumn{1}{c|}{BMFSM} & \multicolumn{1}{c|}{NM} & \multicolumn{1}{c|}{FDWFM} & \multicolumn{1}{c|}{} \\ \hline
\multicolumn{1}{|c|}{(10,10,2,-1)} & \multicolumn{1}{c|}{10} & \multicolumn{1}{l|}{8} & \multicolumn{1}{c|}{7} & \multicolumn{1}{c|}{1.528} & \multicolumn{1}{c|}{1.789} & \multicolumn{1}{c|}{2.534} & \multicolumn{1}{c|}{\begin{tabular}[c]{@{}c@{}}1.0000,1.0000,\\ 0.57735,-0.57735\end{tabular}} \\ \hline
\multicolumn{1}{|c|}{(9.449645,8.198130,1.958279,-2.2299584)} & \multicolumn{1}{c|}{11} & \multicolumn{1}{l|}{8} & \multicolumn{1}{c|}{7} & \multicolumn{1}{c|}{1.687} & \multicolumn{1}{c|}{2.187} & \multicolumn{1}{c|}{2.574} & \multicolumn{1}{c|}{\begin{tabular}[c]{@{}c@{}}1.0000,1.0000,\\ 0.57735,-0.57735\end{tabular}} \\ \hline
\multicolumn{1}{|c|}{(10,10,-1,2)} & \multicolumn{1}{c|}{10} & \multicolumn{1}{l|}{8} & \multicolumn{1}{c|}{7} & \multicolumn{1}{c|}{1.432} & \multicolumn{1}{c|}{1.789} & \multicolumn{1}{c|}{2.341} & \multicolumn{1}{c|}{\begin{tabular}[c]{@{}c@{}}1.0000,1.0000,\\ -0.57735,0.57735\end{tabular}} \\ \hline
\multicolumn{3}{l}{FDWFM-Finite Difference Weerakoon-Fernando Method} & \multicolumn{1}{l}{} & \multicolumn{4}{l}{NM-Newton's method} \\
\multicolumn{3}{l}{BMFSM-Broyden's method for secant method} & \multicolumn{1}{l}{} & \multicolumn{4}{l}{i-Number of iterations to approximate the root} \\
\multicolumn{3}{l}{COC-Computational Order of Convergence} & \multicolumn{1}{l}{} & \multicolumn{1}{l}{} & \multicolumn{1}{l}{} & \multicolumn{1}{l}{} & \multicolumn{1}{l}{}
\end{tabular}
\end{sidewaystable}

\begin{sidewaystable}[]
\centering
\caption{Comparison of FDWFM with Newton\textquotesingle s method and Broyden\textquotesingle s method for systems of six nonlinear equations $f_{1}=x_{1}^{2}+x_{3}^{2}-1, f_{2}=x_{2}^{2}+x_{4}^{2}-1,f_{3}=x_{5}x_{3}^{3}+x_{6}x_{4}^{3}, f_{4}=x_{5}x_{1}^{3}+x_{6}x_{2}^{3}, f_{5}=x_{5}x_{1}x_{3}^{2}+x_{6}x_{4}^{2}x_{2},f_{6}=x_{5}x_{1}^{2}x_{3}+x_{6}x_{2}^{2}x_{4}$}
\label{my-label}
\begin{tabular}{cclccccc}
\hline
\multicolumn{1}{|c|}{\multirow{2}{*}{\begin{tabular}[c]{@{}c@{}}Initial guess\\ X0\end{tabular}}} & \multicolumn{3}{c|}{i} & \multicolumn{3}{c|}{COC} & \multicolumn{1}{c|}{\multirow{2}{*}{Root}} \\ \cline{2-7}
\multicolumn{1}{|c|}{} & \multicolumn{1}{c|}{BMFSM} & \multicolumn{1}{l|}{NM} & \multicolumn{1}{c|}{FDWFM} & \multicolumn{1}{c|}{BMFSM} & \multicolumn{1}{c|}{NM} & \multicolumn{1}{c|}{FDWFM} & \multicolumn{1}{c|}{} \\ \hline
\multicolumn{1}{|c|}{(3,5,4,6,5.5,2,1,-4)} & \multicolumn{1}{c|}{9} & \multicolumn{1}{l|}{7} & \multicolumn{1}{c|}{6} & \multicolumn{1}{c|}{1.632} & \multicolumn{1}{c|}{2.277} & \multicolumn{1}{c|}{2.435} & \multicolumn{1}{c|}{\begin{tabular}[c]{@{}c@{}}0.5368,0.9170,0.8436,\\ 0.3987,0.0000,-0.0000\end{tabular}} \\ \hline
\multicolumn{1}{|c|}{\begin{tabular}[c]{@{}c@{}}(2.5257,5.0538,5.8289,\\ 2.1629,2.4797,-4.9408)\end{tabular}} & \multicolumn{1}{c|}{9} & \multicolumn{1}{l|}{7} & \multicolumn{1}{c|}{6} & \multicolumn{1}{c|}{1.714} & \multicolumn{1}{c|}{2.321} & \multicolumn{1}{c|}{2.498} & \multicolumn{1}{c|}{\begin{tabular}[c]{@{}c@{}}0.5039,0.8519,0.8637,\\ 0.5236,-0.0000,0.0000\end{tabular}} \\ \hline
\multicolumn{1}{|c|}{\begin{tabular}[c]{@{}c@{}}(2.4711,4.3696,6.2511,\\ 1.4369,1,9453,-4.4211)\end{tabular}} & \multicolumn{1}{c|}{9} & \multicolumn{1}{l|}{7} & \multicolumn{1}{c|}{6} & \multicolumn{1}{c|}{1.786} & \multicolumn{1}{c|}{2.625} & \multicolumn{1}{c|}{2.542} & \multicolumn{1}{c|}{\begin{tabular}[c]{@{}c@{}}0.3676,0.9499,0.9299,\\ 0.3123,0.0000,0.0000\end{tabular}} \\ \hline
\multicolumn{4}{l}{FDWFM-Finite Difference Weerakoon-Fernando Method} & \multicolumn{4}{l}{NM-Newton's method} \\
\multicolumn{4}{l}{BMFSM-Broyden's method for secant method} & \multicolumn{4}{l}{i-Number of iterations to approximate the root} \\
\multicolumn{4}{l}{COC-Computational Order of Convergence} & \multicolumn{1}{l}{} & \multicolumn{1}{l}{} & \multicolumn{1}{l}{} & \multicolumn{1}{l}{}
\end{tabular}
\end{sidewaystable}

\begin{sidewaystable}[]
\centering
\caption{Comparison of FDWFM with Newton\textquotesingle s method and Broyden\textquotesingle s method for systems of ten nonlinear equations $f_{1}=x_{1}-0.25428722-0.18324757x_{4}x_{3}x_{9}, f_{2}=x_{2}-0.37842197-0.16275449x_{1}x_{10}x_{6}, f_{3}=x_{3}-0.27162577-0.16955071x_{1}x_{2}x_{10}, f_{4}=x_{4}-0.19807914-0.15585316x_{7}x_{1}x_{6}, f_{5}=x_{5}-0.44166728-1.9950920x_{7}x_{6}x_{3}, f_{6}=x_{6}-0.146541113-0.18922793x_{8}x_{5}x_{10}, f_{7}=x_{7}-0.42937161-0.21180486x_{2}x_{5}x_{8}, f_{8}=x_{8}-0.07056438-0.17081208x_{1}x_{7}x_{6}, f_{9}=x_{9}-0.34504906-0.196127x_{10}x_{6}x_{8}, f_{10}=x_{10}-0.42651102-0.21466544x_{4}x_{8}x_{1}$}
\label{my-label}
\begin{tabular}{cclccccc}
\hline
\multicolumn{1}{|c|}{\multirow{2}{*}{\begin{tabular}[c]{@{}c@{}}Initial guess\\ X0\end{tabular}}} & \multicolumn{3}{c|}{i} & \multicolumn{3}{c|}{COC} & \multicolumn{1}{c|}{\multirow{2}{*}{Root}} \\ \cline{2-7}
\multicolumn{1}{|c|}{} & \multicolumn{1}{c|}{BMFSM} & \multicolumn{1}{l|}{NM} & \multicolumn{1}{c|}{FDWFM} & \multicolumn{1}{c|}{BMFSM} & \multicolumn{1}{c|}{NM} & \multicolumn{1}{c|}{FDWFM} & \multicolumn{1}{c|}{} \\ \hline
\multicolumn{1}{|c|}{(1,1,1,1,1,1,1,1,1,1)} & \multicolumn{1}{c|}{8} & \multicolumn{1}{l|}{5} & \multicolumn{1}{c|}{4} & \multicolumn{1}{c|}{1.234} & \multicolumn{1}{c|}{1.799} & \multicolumn{1}{c|}{2.365} & \multicolumn{1}{c|}{\begin{tabular}[c]{@{}c@{}}0.2578,0.3810,0.2878,0.2006\\ 0.4452,0.1491,0.4320,0.0734\\ 0.3459,0.4273\end{tabular}} \\ \hline
\multicolumn{1}{|c|}{\begin{tabular}[c]{@{}c@{}}(-0.3956,-1.3108,-0.3927,4.8163,\\ -3.4359,3.555,1.4476,-1.2372,\\ -3.0907,-0.7174)\end{tabular}} & \multicolumn{1}{c|}{12} & \multicolumn{1}{l|}{10} & \multicolumn{1}{c|}{8} & \multicolumn{1}{c|}{1.435} & \multicolumn{1}{c|}{1.827} & \multicolumn{1}{c|}{2.481} & \multicolumn{1}{c|}{\begin{tabular}[c]{@{}c@{}}0.3452,1.2453,1,0342,1,4352,\\ 1.4533,1.4563,1.8453,1.7453,\\ 1.3435,1.8464\end{tabular}} \\ \hline
\multicolumn{1}{|c|}{\begin{tabular}[c]{@{}c@{}}(1.625,1.780,1.0811,1.9293,\\ 1.7757,1.4867,1.4358,1.4467\\ 1.3063,1.5085)\end{tabular}} & \multicolumn{1}{c|}{10} & \multicolumn{1}{l|}{8} & \multicolumn{1}{c|}{6} & \multicolumn{1}{c|}{1.654} & \multicolumn{1}{c|}{1.958} & \multicolumn{1}{c|}{2.499} & \multicolumn{1}{c|}{\begin{tabular}[c]{@{}c@{}}1.8430,1.9683,1.6191,2.0850,\\ 2.5636,2.4194,2.7151,2.1386,\\ 2.5682,2.1907\end{tabular}} \\ \hline
\multicolumn{4}{l}{FDWFM-Finite Difference Weerakoon-Fernando Method} & \multicolumn{4}{l}{NM-Newton's method} \\
\multicolumn{4}{l}{BMFSM-Broyden's method for secant method} & \multicolumn{4}{l}{i-Number of iterations to approximate the root} \\
\multicolumn{4}{l}{COC-Computational Order of Convergence} & \multicolumn{1}{l}{} & \multicolumn{1}{l}{} & \multicolumn{1}{l}{} & \multicolumn{1}{l}{}
\end{tabular}
\end{sidewaystable}

\newpage
Table 1 shows the result for comparison of FDWFM with secant method and Newton's for uni-variate nonlinear equations with real roots. Here computer order of convergence for FDWFM is near 2.4. It is higher than Secant method and Newton method. Then we applied the method for complex roots.
When the method was applied to nonlinear equations for complex roots it also gives the same satisfactory results. Then the method was extended to systems of nonlinear equations.
When we apply the method to systems of non-linear equations there were some difficulties. In Improved Newton’s method it was necessary to obtain Jacobian matrices to solve systems of nonlinear equations but for this resulting method the Jacobian matrix becomes a vector. So the challenge was to find the inverse of the vector. So Broyden’s method was used to overcome that difficulty. After following the technique used in Broyden’s method resulting formula for systems of nonlinear equations with two three four six and ten variables, the results were again encouraging just like for the one variable case. So the objective was achieved there as well. However, when we follow Broyden’s method we cannot do away with the derivative part. \\

\section{Conclusion}
For all nonlinear equations we have considered, the new algorithm FDWFM seems to be more efficient than any other algorithm without the derivatives to numerically solve them. This method gives a computational order of convergence higher than any other existing method in the absence of the derivatives. It returns the same order of convergence for nonlinear equations with complex roots as well as for systems of nonlinear equations. Apparently, the FDWFM had three function evaluations and secant method had two function evaluations. But according to the computed results (Table 1 and Table 2) most of the results, total number of function evaluations required is less than or equal that secant method. Thus FDWFM can be considered as a superior method giving faster convergence to find the roots of nonlinear equations in the absence of the derivative for uni-variate nonlinear equations with complex roots as well as for the multivariate systems of nonlinear equations. Since the proposed method is even faster than the universally accepted second order Newton's method while meeting the requirement of not having the derivative of the function as well, this algorithm will undoubtedly be very useful to the scientific and the industrial community.

\section{References}
	
\begin{enumerate}
\item Atkinson K. (1988). \textit{An Introduction to Numerical Analysis}, 2nd edition, University of lowa.
\item Burden R.L. \& Faires J.D. (2005). \textit{Numerical Analysis}, 8th edition, Bob Pirtle, USA.
\item Dennis J.E. \& Robert B.S. (1987) \textit{ Numerical Methods for Unconstrained Optimization and Nonlinear Equation}, Society for Industrial and Applied Mathematics, USA.
\item Fernando T. G. I. (1998). \textit {Improved Newton's Method for Solving Nonlinear Equations}, M. Sc. industrial Mathematics Thesis, University of Sri Jayewardenepura. \\  DOI: https://www.researchgate.net
\item Nishani H. P. S., Weerakoon S., Fernando T.G.I. \& Liyanage M. (2014). \textit {Third order convergence of Improved Newton's method for systems of nonlinear equations}, 502/E1, Proceedings of the annual sessions of Sri Lanka Association for the Advancement of Science, Sri Lanka.
\item	Nishani H.P.S. (2015). \textit{A Variant of Newton\textquotesingle s Method with Accelerated Third Order Convergence for Systems of Nonlinear Equations}, M. Sc in Industrial Mathematics thesis, University of Sri Jayewardenepura.

\item	Nishani H.P.S., Weerakoon S., Fernando T.G.I. \& Liyanage L.M. (2018). \textit{Weerakoon Fernando Method with accelerated third order convergence for systems of nonlinear equations}, International Journal of Mathematical Modelling and Numerical Optimisation \textbf{8}(3):287 - 304\\
DOI:www.inderscience.com/info/inarticle.php?artid=89010

\item Weerakoon S. \& Fernando T.G.I. (2000). \textit{A Variant of Newton\textquotesingle s Method with Accelerated Third-Order Convergence}, Applied Mathematics Letters \textbf{13} (8)\\  DOI: http://www.sciencedirect.com/science/article/pii5  
\item Young T. \& Martin J. (2015). \textit{Introduction to Numerical Methods and Matlab
Programming for Engineers}, Department of Mathematics, University of Ohio.

\end{enumerate}

\end{document}